# The Göbbellian Syndrome

## Can we *really* falsify truth by dictat?

Bhupinder Singh Anand[1]

## 1. Introduction – Non-standard models of PA

In a 1996 talk on the Gödelian argument [Lu96], J.R. Lucas commented as follows:

> … in the case of First-order Peano Arithmetic there are Gödelian formulae (many, in fact infinitely many, one for each system of coding) which are not assigned truth-values by the rules of the system, and which could therefore be assigned either TRUE or FALSE, each such assignment yielding a logically possible, consistent system. These systems are random vaunts, all satisfying the core description of Peano Arithmetic.

## 2. A Göbbellian dictat?

Can we *really* falsify truth by such a Göbbellian dictat?

In other words, if $[(Ax)R(x)]$ is the PA-unprovable Gödelian formula[2] - which is true under the standard interpretation - can we add its negation, $[\sim(Ax)R(x)]$, as an axiom to PA, and still obtain a consistent system?

In his seminal 1931 paper, Gödel argued ([Go31], p27) that, if an arithmetic such as PA is $\omega$-consistent, then the system PA+$[\sim(Ax)R(x)]$, say PA*, is consistent, but not $\omega$- consistent. He asserted, without proof, that since $[(Ax)R(x)]$ is unprovable in his arithmetic, P, it continues to remain unprovable in $P+[\sim(Ax)R(x)]$. This, of course, would not be the case if the axiomatic addition of $[\sim(Ax)R(x)]$ to P were to introduce an inconsistency.

---

[1] The author is an independent scholar. E-mail: re@alixcomsi.com; anandb@vsnl.com. Postal address: 32, Agarwal House, D Road, Churchgate, Mumbai - 400 020, INDIA. Tel: +91 (22) 2281 3353. Fax: +91 (22) 2209 5091.

[2] This is the formula that Gödel defines, and refers to, by its Gödel-number, 17*Genr* ([Go31] p25).

## 3. Does ZF model PA*?

Classically, any first order theory with equality is consistent if, and only if, it has a model ([Me64]; p51(II) & p65, proposition 2.12). Gödel's postulation of the relative consistency of PA*, therefore, implies that PA* has a model.

Further, an implicit belief of classical theory is that first order mathematical theories can be interpreted suitably in ZF, and so every model of ZF is a model of such systems.

The question, thus, arises: Is PA* such a system?

Now, standard interpretations of classical theory argue ([Me64], p192-193) that the axioms of PA can be interpreted in ZF - in the sense that, under suitable interpretations of the primitive symbols of PA in ZF, the axioms of PA can be transformed into theorems of ZF, and the rules of inference preserved, by restricting the variables to Cantor's first transfinite ordinal, $\omega$, so that the PA-formula $[(Ax)R(x)]$, for instance, would translate as $[(Ax)((x \in \omega) => R(x))]$.

Clearly, every model of ZF is, thus, a model of the interpreted axioms of PA. However, does this assure us of a ZF model for PA*?

Now, for PA* to have a model in ZF, we would need $[\sim(Ax)((x \in \omega) => R(x))]$ to be a theorem of ZF. This, of course, is not possible if ZF is consistent; Gödel has shown that $[(Ax)((x \in \omega) => R(x))]$ is a theorem of ZF, since it is true in every model of ZF.

## 4. Does PA* have a model outside ZF?

Further, the above argument would hold for every model of PA*, since the induction axiom of PA would hold in a non-standard model (cf. [Me64], p115, Ex.2) if, and only if - as in ZF - we are able to restrict the range of the variables in it to the natural numbers.

So, was Gödel's postulation simply a case of a wrongly proven conjecture, or are there alternative arguments for concluding that PA* is consistent, and that it has a non-trivial non-standard model?

## 5. In conclusion

We note that alternative non-standard models of PA, such as those in Mendelson ([Me64], p116, Ex.2), or even Cantor's ordinal numbers, 0, 0', 0", ..., ω, ω', ω", ..., ω.0', (ω.0")', (ω.0")", ..., essentially consist of a sequence of equivalence classes, {0, ω, ω.0", ...}, {0', ω', (ω.0")', ...}, {0", ω", (ω.0")", ...}, ...; such sequences are, themselves, simply contractions ([Me64], p80) - of non-standard models of PA - that are isomorphic to the standard model of PA.

This follows since each of {ω, ω', ω", ...}, {ω.0",( ω.0")', (ω.0")", ...}, etc. is isomorphic to the standard model, {0, 0', 0", ...}, of PA.

It is difficult to see the significance given to such, non-standard, models in standard interpretations of classical theory, since any model of any theory can be used to create non-standard models - of the above kind - by simply re-naming the primitive constants of the theory, and arranging the renamed domains sequentially.

Clearly, such extended domains are trivial, since they cannot have any mathematical properties - pertaining to the theory of which they are the postulated models - that are not, already, in the standard model of the theory.

## References

[Go31]  Gödel, Kurt. 1931. *On formally undecidable propositions of Principia Mathematica and related systems I.* Translated by Elliott Mendelson. In M. Davis (ed.). 1965. The Undecidable. Raven Press, New York.

[Lu61]  Lucas, John, R. 1961. *Minds, Machines and Gödel.* First published in *Philosophy*, XXXVI, 1961, pp.(112)-(127); reprinted in *The Modeling of Mind*, Kenneth M. Sayre and Frederick J. Crosson, eds., Notre Dame Press, 1963, pp.[269]-[270]; and *Minds and Machines*, ed. Alan Ross Anderson, Prentice-Hall, 1954, pp.{43}-{59}.

<*Web page*: http://users.ox.ac.uk/~jrlucas/Godel/mmg.html>

[Lu96]  Lucas, John, R. 1996. *The Gödelian Argument: Turn Over the Page*. A talk given on 25/5/96 at a BSPS conference in Oxford.

<*Web page*: http://users.ox.ac.uk/~jrlucas/Godel/turn.html>

[Me64]  Mendelson, Elliott. 1964. Introduction to Mathematical Logic. Van Norstrand, Princeton.
(*Updated: Wednesday*, 15$^{th}$ *June 2005, 2:23:58 AM IST*, *by re@alixcomsi.com*)